\documentclass{amsart}
\usepackage{graphicx}
\usepackage{amsfonts}

\begin{document}

\title[THE NUMERICAL GENERALIZED LEAST-SQUARES ESTIMATOR OF ...]
{THE NUMERICAL GENERALIZED LEAST-SQUARES ESTIMATOR OF
AN UNKNOWN CONSTANT MEAN OF 
RANDOM FIELD}

\author{T. SUS{\L}O}
\email{tomasz.suslo@gmail.com}

\begin{abstract}
We constraint on computer the best linear unbiased generalized statistics 
of random field for the best linear unbiased generalized statistics of 
an unknown constant mean of random field and derive the numerical generalized least-squares estimator of an unknown constant mean of random field.
We derive the third constraint of spatial statistics and show 
that the classic generalized least-squares estimator of an unknown
constant mean of the field is only an asymptotic disjunction of 
the numerical one.
\end{abstract}

\maketitle

\thispagestyle{empty}

\section{The best linear unbiased generalized statistics}
\noindent
{\bf Remark.}
To simplify notation we use Einstein summation convention then
$$
\sum_{i=1}^n \omega^i_j \rho_{ij}
=
\omega^i_j \rho_{ij}
=
w' r \,
$$
where
$$
\begin{array}{cccccccccc}
w
&
=
&
\underbrace{
\left[
\begin{array}{c}
\omega_j^{\it 1} \\
\vdots \\
\omega_j^{\it n} \\
\end{array}
\right]
}_{n \times 1} \ ,
&
r
&
=
&
\underbrace{
\left[
\begin{array}{c}
\rho_{\it 1j} \\
\vdots \\
\rho_{\it nj} \\
\end{array}
\right]
}_{n\times 1}
\end{array}
$$
are given vectors and
$$
\sum_{i=1}^n\omega^i_j\sum_{l=1}^n \rho_{il} \omega^l_j
=
\omega^i_j \rho_{il} \omega^l_j
=
w' \Lambda w \ ,
$$
where
$$
\begin{array}{ccccc}
\Lambda
&
=
&
{\underbrace{
\left[
\begin{array}{ccc}
\rho_{\it 11} & \ldots & \rho_{\it 1n} \\
\vdots & \ddots & \vdots \\
\rho_{\it n1}& \ldots & \rho_{\it nn} \\
\end{array}
\right]}_{n \times n}}
\end{array}
$$
is given matrix.

Let us consider the random field $V_j;~j \in {\mathbb N}_1$
with an unknown constant mean $m$ and variance $\sigma^2$
its estimation statistics $\hat{V}_j$
and the variance of the difference $R_j=V_j-\hat{V}_j$,
where $E\{V_j\}=E\{\hat{V}_j\}=m$,
as covariance
$$
\begin{array}{ccc}
D^2\{V_j-\hat{V}_j\}
&
=
&
Cov\{(V_j-\hat{V}_j)(V_j-\hat{V}_j)\}
\\

&
=
&
Cov\{V_jV_j\}
-
Cov\{V_j\hat{V}_j\}
-
Cov\{\hat{V}_jV_j\}
+
Cov\{\hat{V}_j\hat{V}_j\}
\\

&
=
&
Cov\{V_jV_j\}
-
2Cov\{\hat{V}_jV_j\}
+
Cov\{\hat{V}_j\hat{V}_j\}
\end{array}
$$
and the linear estimation statistics (weighted variable)
$\hat{V}_j
=
\sum_{i=1}^n \omega^i_j V_i
=
\omega^i_j V_i;~j \subset i=1,\ldots,n$
at $j\ge n+1$ then
\begin{eqnarray}
D^2\{R_j\}
&
=
&
Cov\{V_jV_j\}-2Cov\{\hat{V}_jV_j\}+Cov\{\hat{V}_j\hat{V}_j\} \nonumber
\\
&
=
&
Var\{V_j\}
-2Cov\{\sum_i\omega_j^i V_i V_j\}
+
Cov\{(\sum_i\omega^i_j V_i)(\sum_i\omega^i_j V_i)\} \nonumber
\\
&
=
&
\sigma^2-2\sum_i\omega^i_j Cov\{V_i V_j\}
+
\sum_i\sum_l\omega^i_j \omega^l_j Cov\{V_i V_l\} \nonumber
\\
&
=
&
\sigma^2
-
2\sigma^2|\omega^i_j \rho_{ij}|
+
\sigma^2 |\omega^i_j \rho_{il} \omega^l_j| \nonumber
\\
&
=
&
\sigma^2
\pm 2\sigma^2\omega^i_j \rho_{ij}
\mp \sigma^2 \omega^i_j \rho_{il} \omega^l_j \ ,
\label{D}
\end{eqnarray}
where $\rho_{ij};~i=1,\ldots,n$ is given
vector of correlations and $\rho_{il};~i,l=1,\ldots,n$
is given (symmetric) matrix of correlations~(see~Appendix~\ref{sec:A}). \\
The unbiasedness constraint
(the first constraint on the estimation statistics)
$$
E\{R_j\}=E\{V_j-\hat{V}_j\}
=
E\{V_j\}-E\{\hat{V}_j\}
=
E\{V_j\}-E\{\omega^i_j V_i\}
=
m-m\sum_{i=1}^n \omega^i_j
=
0
$$
equal to 
\begin{equation}
\sum_{i=1}^n \omega^i_j=
f_{\it 1i} \omega^i_j
=
\omega^i_j f_{\it i1}
=1
\label{uc}
\end{equation}
gives the first equation
$$
\begin{array}{ccccccccc}
{\underbrace{
\left[
\begin{array}{ccc}
1 & \ldots & 1
\end{array}
\right]}_{1 \times n}}
&
\cdot
&
\underbrace{
\left[
\begin{array}{c}
\omega^1_j \\
\vdots \\
\omega^n_j \\
\end{array}
\right]
}_{n \times 1}
&
=
&
{\underbrace{
\left[
\begin{array}{ccc}
\omega^1_j & \ldots & \omega^n_j
\end{array}
\right]}_{1 \times n}}
&
\cdot
&
\underbrace{
\left[
\begin{array}{c}
1 \\
\vdots \\
1 \\
\end{array}
\right]
}_{n \times 1}
&
=
&
1 \ .
\end{array}
$$
The minimization constraint
(the second constraint on the estimation statistics
-- the statistics is the best)
\begin{equation}
\frac{\partial D^2\{R_j\}}{\partial \omega^i_j}
=
\pm 2\sigma^2\rho_{ij}
\mp 2\sigma^2\rho_{il}\omega^l_j
\mp 2\sigma^2f_{\it i1}\mu_j^{\it 1}
= 0 \ ,
\label{mc}
\end{equation}
where~(\ref{D})
$$
D^2\{R_j\}=\sigma^2
\pm 2\sigma^2\omega^i_j\rho_{ij}
\mp \sigma^2\omega^i_j\rho_{il} \omega^l_j
\mp 2\sigma^2\underbrace{\left(\omega^i_j f_{\it i1}
- 1\right)}_0 \mu_j^{\it 1} \ ,
$$
produces $n$ equations in $n+1$ unknowns
the kriging weights $\omega_j^i$ and a Lagrange
parameter $\mu^{\it 1}_j$
$$
\begin{array}{cccccl}
{\underbrace{
\left[
\begin{array}{cccc}
\rho_{\it 11} & \ldots & \rho_{\it 1n} & 1 \\
\vdots & \ddots & \vdots & \vdots \\
\rho_{\it n1} & \ldots & \rho_{\it nn} & 1 \\
\end{array}
\right]}_{n\times(n+1)}}
&
\cdot
&
\underbrace{
\left[
\begin{array}{c}
\omega_j^{\it 1} \\
\vdots \\
\omega_j^{\it n} \\
\mu_j^{\it 1} \\
\end{array}
\right]
}_{(n+1) \times 1}
&
=
&
\underbrace{
\left[
\begin{array}{c}
\rho_{\it 1j} \\
\vdots \\
\rho_{\it nj} \\
\end{array}
\right]
}_{n \times 1}
\end{array}
$$
this system of equations if multiplied by $\omega^i_j$
$$
\omega^i_j \rho_{il} \omega^l_j
+
\underbrace{\omega^i_j f_{\it i1}}_1 \mu^{\it 1}_j
=
\omega^i_j \rho_{ij} \ ,
$$
and substituted into
$$
\begin{array}{ccc}
D^2\{R_j\}
&
=
&
E\{[V_j-\hat{V}_j]^2\}-\underbrace{E^2\{V_j-\hat{V}_j\}}_0 \\

&
=
&
E\{[(V_j-m)-(\hat{V}_j-m)]^2\} \\

&
=
&
E\{[V_j-m]^2\}-2(E\{V_j\hat{V}_j\}-m^2)+E\{[\hat{V}_j-m]^2\} \\

&
=
&
\sigma^2
-2 \sigma^2|\omega^i_j\rho_{ij}|
+ \sigma^2|\omega^i_j\rho_{il} \omega^l_j| \\

&
=
&
\sigma^2
\pm 2 \sigma^2 \omega^i_j \rho_{ij}
\mp \sigma^2 \omega^i_j \rho_{il} \omega^l_j 
\end{array}
$$
since variance of the (estimation) statistics is minimized
\begin{eqnarray}
E\{[\hat{V}_j-m]^2\}
&
=
&
Cov\{(\omega^i_j V_i)(\omega^i_j V_i)\} \nonumber
\\
&
=
&
\sum_i\sum_l\omega^i_j\omega^l_j Cov\{V_i V_l\} \nonumber
\\
&
=
&
\sigma^2|\omega^i_j \rho_{il} \omega^l_j| \nonumber
\\
&
=
&
\mp \sigma^2 \omega^i_j \rho_{il} \omega^l_j \nonumber
\\
&
=
&
\mp \sigma^2 (\omega^i_j \rho_{ij}-\mu_j^{\it 1})
\label{o}
\end{eqnarray}
gives
\begin{equation}
D^2\{R_j\}
=
E\{[V_j-\hat{V}_j]^2\}
=
E\{[(V_j-m)-(\hat{V}_j-m)]^2\}
=
\sigma^2
(1 \pm (\omega^i_j \rho_{ij}+\mu^{\it 1}_j))
\label{oo}
\end{equation}
the constraints~(\ref{uc}) and~(\ref{mc}) produce $n+1$ equations
in $n+1$ unknowns
\begin{equation}
\begin{array}{cccccl}
{\underbrace{
\left[
\begin{array}{cccc}
\rho_{\it 11} & \ldots & \rho_{\it 1n} & 1 \\
\vdots & \ddots & \vdots & \vdots \\
\rho_{\it n1} & \ldots & \rho_{\it nn} & 1 \\
1 & \ldots & 1 & 0 \\
\end{array}
\right]}_{(n+1) \times (n+1)}}
&
\cdot
&
\underbrace{
\left[
\begin{array}{c}
\omega_j^{\it 1} \\
\vdots \\
\omega_j^{\it n} \\
\mu^{\it 1}_j \\
\end{array}
\right]
}_{(n+1) \times 1}
&
=
&
\underbrace{
\left[
\begin{array}{c}
\rho_{\it 1j} \\
\vdots \\
\rho_{\it nj} \\
1 \\
\end{array}
\right]
}_{(n+1) \times 1} \ .
\end{array}
\label{ke}
\end{equation}

\section{The classic best linear unbiased generalized statistics of
an unknown constant mean of the field}
\noindent
{\bf Remark.} When we consider an independent set of
the random variables $V_i;~i=1,\ldots,n$
with an unknown constant mean $m$ and variance $\sigma^2$
the best linear unbiased ordinary (estimation) statistics
$\hat{V}_j=\omega^i_j V_i$ of the field $V_j;~j \subset i=1,\ldots,n$
has the asymptotic property
\begin{equation}
\lim_{n \rightarrow \infty}
E\{[\omega^i_j V_i-m]^2\}
=
0
\label{sal}
\end{equation}
whilst for spatial dependence between random variables
(the best linear unbiased generalized statistics)
we get~(see~Appendix~\ref{sec:B})
\begin{equation}
\lim_{n \rightarrow \infty}
\lim_{j \rightarrow \infty}
E\{[\omega^i_j V_i-m]^2\}
=
0 \ .
\label{tal}
\end{equation}
Due to different asymptotic limits between~(\ref{sal}) and~(\ref{tal})
the ordinary least-squares estimator of an unknown constant mean $m$
of the field, the best linear unbiased estimator of
an unknown constant mean $m$ of the field,
can not be so easy generalized (like it was in past).

Let us constraint the best linear unbiased
generalized (estimation) statistics $\hat{V}_j=\omega^i_j V_i$
of the random field $V_j;~j \subset i=1,\ldots,n$,
when for finite $n$ and $j \rightarrow \infty$
the vector of correlations simplifies to
\begin{equation}
\underbrace{
\left[
\begin{array}{c}
\rho_{1j} \\
\vdots \\
\rho_{nj} \\
\end{array}
\right]
}_{n \times 1}
=
\xi
\underbrace{
\left[
\begin{array}{c}
1 \\
\vdots \\
1 \\
\end{array}
\right]
}_{n \times 1}
\qquad \xi \rightarrow 0^-~(\mbox{or} ~\xi \rightarrow 0^+)
\label{cv}
\end{equation}
then from~(\ref{uc})
\begin{equation}
\lim_{j \rightarrow \infty}
\omega^i_j \rho_{ij}=\xi \omega^i_j f_{\it i1} = \xi
\label{wrho}
\end{equation}
it holds~(\ref{oo})
\begin{equation}
\lim_{j \rightarrow \infty}
E\{[V_j-\omega^i_j V_i]^2\}
=
\lim_{j \rightarrow \infty}
\sigma^2
(1 \pm (\omega^i_j \rho_{ij}+\mu^{\it 1}_j))
=
\sigma^2
(1 \pm (\xi + \mu^{\it 1}_j))
\label{joo}
\end{equation}
for the co-ordinate independent statistics
of an unknown constant mean of the field $V_j$
with the constraint on~(\ref{joo})
\begin{equation}
\lim_{j \rightarrow \infty}
E\{[V_j-\omega^i_j V_i]^2\}
=
\sigma^2
=
E\{[V_j-m]^2\}
\label{tc}
\end{equation}
given by constrained
from~(\ref{joo})
\begin{equation}
\mu_j^{\it 1}=-\xi
\label{mu}
\end{equation}
and from~(\ref{cv}) the system of equations~(\ref{ke})
$$
\begin{array}{cccccl}
{\underbrace{
\left[
\begin{array}{cccc}
\rho_{\it 11} & \ldots & \rho_{\it 1n} & 1 \\
\vdots & \ddots & \vdots & \vdots \\
\rho_{\it n1} & \ldots & \rho_{\it nn} & 1 \\
1 & \ldots & 1 & 0 \\
\end{array}
\right]}_{(n+1) \times (n+1)}}
&
\cdot
&
\underbrace{
\left[
\begin{array}{c}
\omega_j^{\it 1} \\
\vdots \\
\omega_j^{\it n} \\
- \xi \\
\end{array}
\right]
}_{(n+1) \times 1}
&
=
&
\underbrace{
\left[
\begin{array}{c}
\xi \\
\vdots \\
\xi \\
1 \\
\end{array}
\right]
}_{(n+1) \times 1}
&
\end{array}
$$
equivalent to
$$
\Lambda w -\xi F = \xi F
$$
and
$$
F'w=1 \ ,
$$
where
$$
\begin{array}{cccccccccccc}
w
&
=
&
\underbrace{
\left[
\begin{array}{c}
\omega_j^{\it 1} \\
\vdots \\
\omega_j^{\it n} \\
\end{array}
\right]
}_{n \times 1} \ ,
&
F
&
=
&
\underbrace{
\left[
\begin{array}{c}
1 \\
\vdots \\
1 \\
\end{array}
\right]
}_{n \times 1}
&

,
&
\Lambda
&
=
&
\Lambda'
&
=
&
{\underbrace{
\left[
\begin{array}{ccc}
\rho_{\it 11} & \ldots & \rho_{\it 1n} \\
\vdots & \ddots & \vdots \\
\rho_{\it n1} & \ldots & \rho_{\it nn} \\
\end{array}
\right]}_{n \times n}} \ ,
\end{array}
$$
with the solution
\begin{equation}
\xi=\frac{1}{2 F' \Lambda^{-1} F}
\label{xii}
\end{equation}
and
\begin{equation}
w=\frac{\Lambda^{-1}F}{F' \Lambda^{-1} F}
\end{equation}
of the classic best linear unbiased generalized statistics
for finite $n$ and $j \rightarrow \infty$ of
an unknown constant mean of the field
\begin{equation}
\lim_{j \rightarrow \infty}
w' V=\frac{F'\Lambda^{-1}V}{F'\Lambda^{-1} F} \ ,
\label{csol}
\end{equation}
where
$$
V
=
\underbrace{
\left[
\begin{array}{c}
V_1 \\
\vdots \\
V_n \\
\end{array}
\right]
}_{n \times 1} \ ,
$$
with constrained minimized variance of
the best linear unbiased generalized (estimation) statistics~(\ref{o})
as its variance~(from(\ref{wrho})~and~(\ref{mu}))
$$
\lim_{j \rightarrow \infty}
E\{[\omega^i_j V_i-m]^2\}
=
\lim_{j \rightarrow \infty}
\mp\sigma^2(\omega^i_j \rho_{ij}-\mu_j^{\it 1})
=
\mp\sigma^2 (\xi - \mu_j^{\it 1})
=
\mp\sigma^2 2\xi
$$
then~(from(\ref{xii}))
$$
\lim_{j \rightarrow \infty}
E\{[w'V-m]^2\}
=
\mp\sigma^2 2\xi
=
\frac{\mp\sigma^2}{F'\Lambda^{-1}F} \ ,
$$
with the classic generalized least-squares estimator
for finite $n$ and $j \rightarrow \infty$ of an unknown
constant mean $m$ of the field
\begin{equation}
\lim_{j \rightarrow \infty}
w' {\bf v}=\frac{F' \Lambda^{-1}{\bf v}}{F'\Lambda^{-1}F}
\label{cglse}
\end{equation}
based on observation ${\bf v}$ seen as outcome of $V$.

\section{The numerical best linear unbiased generalized statistics of
an unknown constant mean of the field}
\begin{figure}
\includegraphics[width=12cm,height=6cm]{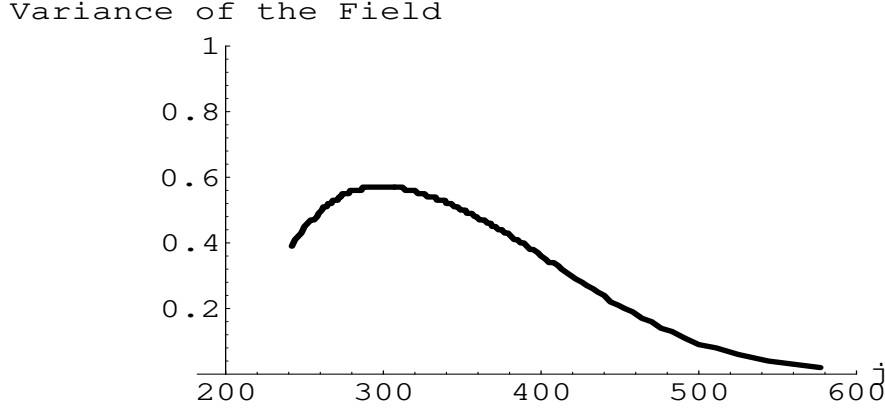}
\caption{\label{Fig1} Variance of 
the numerical best linear unbiased generalized statistics
for finite $n$ at finite $j \ge n+1=182+1$ of
an unknown constant mean $m=E\{V_j\}$ of
the field $V_j$ in units of the variance $\sigma^2$ of
the field computed 139 times for 
the negative correlation function~(\ref{cf}) with
the parameter $t=182+1,\ldots,182+139$.}
\end{figure}

\noindent
To remove the asymptotic limit of
the classic best linear unbiased generalized statistics
for finite $n$ and $j \rightarrow \infty$ of
an unknown constant mean $m=E\{V_j\}$ of
the field $V_j$ with the constraint~(\ref{tc})
$$
\lim_{j \rightarrow \infty}E\{[V_j-\omega^i_j V_i]^2\}
=
\sigma^2
=
E\{[V_j-m]^2\} \ ,
$$
the best linear unbiased generalized (estimation) statistic of
the field $V_j;~j \subset i=1,\ldots,n$ at finite $j \ge n+1=182+1$
$$
\hat{V}_j = \sum_{i=1}^{n=182} \omega^i_j V_i
$$
given by the kriging algorithm~(\ref{ke}) for $n=182$
$$
\begin{array}{cccccl}
\underbrace{
\left[
\begin{array}{c}
\omega_j^{\it 1} \\
\vdots \\
\omega_j^{\it n} \\
\mu^{\it 1}_j \\
\end{array}
\right]
}_{(n+1) \times 1}
&
=
{\underbrace{
\left[
\begin{array}{cccc}
\rho_{\it 11} & \ldots & \rho_{\it 1n} & 1 \\
\vdots & \ddots & \vdots & \vdots \\
\rho_{\it n1} & \ldots & \rho_{\it nn} & 1 \\
1 & \ldots & 1 & 0 \\
\end{array}
\right]}_{(n+1) \times (n+1)}}^{-1}
&
\cdot
&
\underbrace{
\left[
\begin{array}{c}
\rho_{\it 1j} \\
\vdots \\
\rho_{\it nj} \\
1 \\
\end{array}
\right]
}_{(n+1) \times 1}
\end{array}
$$
the negative correlation function with the parameter
$t=182+1,\ldots,182+139$
\begin{equation}
\rho(\Delta_{ij})=\left\{
     \begin{array}{ll}
     -1 \cdot {\displaystyle t}^{-0.62590\displaystyle [\Delta_{ij} \slash t]^2},&
      \qquad \mbox{for}~~\Delta_{ij}=|i-j|> 0,\\
      +1, & \qquad \mbox{for}~~\Delta_{ij}=|i-j|=0,\\
      \end{array}
      \right.
\label{cf}
\end{equation}
was constrained~(from~(\ref{oo}))
on computer (139 times)
for the numerical best linear unbiased generalized statistics
for finite $n$ at finite $j$ of
an unknown constant mean $m=E\{V_j\}$ of
the field $V_j$ with the third constraint of spatial statistics
\begin{equation}
E\{[V_j-\omega^i_j V_i]^2\}
=
\sigma^2
=
E\{[V_j-m]^2\}
\label{coss}
\end{equation}
equivalent to
\begin{equation}
\omega^i_j \rho_{ij}+ \mu_j^{\it 1}=0
\label{ce}
\end{equation}
with constrained minimized variance of the best linear unbiased generalized 
(estimation) statistics~(\ref{o}) as its variance (see~Fig.~\ref{Fig1}).

Our aim was to derive
for the negative correlation function~(\ref{cf}) with
the parameter $t=182+1,\ldots,182+139$
the numerical generalized least-squares estimator
$\omega^i_j v_i$ of
an unknown constant mean $m=E\{V_j\}$ of the field $V_j$
in fact the proper best linear unbiased (generalized)
estimator of an unknown constant mean $m=E\{V_j\}$ of the field $V_j$
\begin{figure}
\includegraphics[width=12cm,height=6cm]{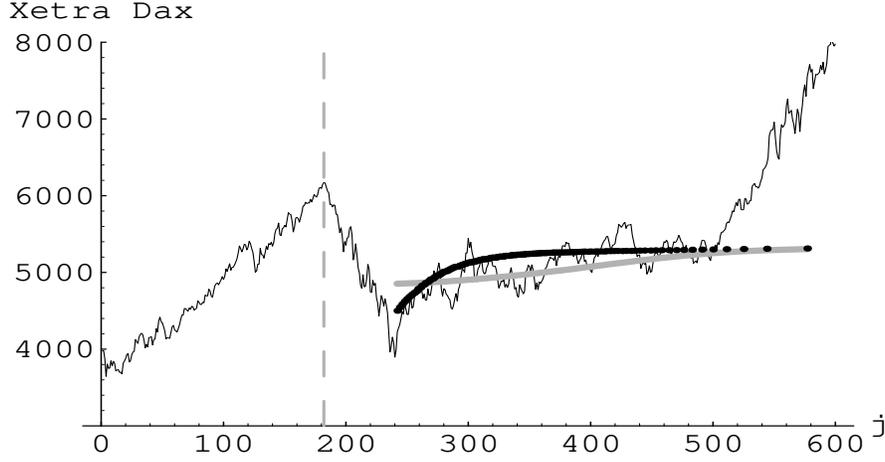}
\caption{\label{Fig2} Long-lived asymmetric index profile,Xetra Dax Index from 23 X 1997 up to 10 III 2000 (600 close quotes)
the numerical generalized least-squares estimator
$\omega^i_j v_i$ of an unknown constant mean $m=E\{V_j\}$ of
the field $V_j;~j \subset i=1,\ldots,182$ (black dots) based on $v_i=v_1,\ldots,v_{182}$ is compared for the negative correlation function~(\ref{cf}) with the parameter $t=182+1,\ldots,182+139$ at 
finite $j \ge n+1=182+1$ to 
the classic generalized least-squares
estimator $\lim_{j \rightarrow \infty}\omega^i_j v_i$ of 
an unknown constant mean $m=E\{V_j\}$ of the field $V_j$
(grey line) with the same correlation function and based on 
the same sample. 
The classic estimator is the first approximation of 
the numerical estimator at final $j=577$ for final $t=182+139$.
The dashed vertical line represents $j=n=182$.}
\end{figure}
given at finite $j \ge n+1=182+1$
by numerical approximation to root of
the equation~(\ref{ce}).
This (co-ordinate dependent) generalized least-squares estimator
$\omega^i_j v_i$ was compared to the (co-ordinate independent)
classic generalized least-squares estimator
$\lim_{j \rightarrow \infty}\omega^i_j v_i$
of an unknown constant mean of the field~(\ref{cglse})
$$
\lim_{j \rightarrow \infty}
w'{\bf v}=\frac{F'\Lambda^{-1}{\bf v}}{F' \Lambda^{-1} F}
$$
based on the same observation
an initial amplification $v_i=v_1,\ldots,v_{182}$ of
long-lived asymmetric index profile recorded
by $600$ close quotes of Xetra Dax Index shown
in~Fig.~\ref{Fig2} then
$$
{\bf v}
=
\underbrace{
\left[
\begin{array}{c}
v_1 \\
\vdots \\
v_{182} \\
\end{array}
\right]
}_{n \times 1}
$$
with the same correlation function~(\ref{cf}).

Since the classic best linear unbiased generalized statistics
for finite $n$ and $j \rightarrow \infty$ of
an unknown constant mean $m=E\{V_j\}$ of the field $V_j$
with the constraint
$$
\lim_{j \rightarrow \infty}
E\{[V_j-\omega^i_j V_i]^2\}
=
\sigma^2
=
E\{[V_j-m]^2\} \ ,
$$
is an asymptotic disjunction for $j \rightarrow \infty$ of
the numerical best linear unbiased generalized statistics for
finite $n$ at finite $j$ of an unknown constant mean $m=E\{V_j\}$ of
the field $V_j$ with the constraint
$$
E\{[V_j-\omega^i_j V_i]^2\}
=
\sigma^2
=
E\{[V_j-m]^2\} \ ,
$$
then the correct classic generalized least-squares estimator
$\lim_{j \rightarrow \infty} \omega^i_j v_i$ of
an unknown constant mean $m$ of the field is an asymptotic
disjunction for $j \rightarrow \infty$
of the numerical generalized least-squares estimator
$\omega^i_j v_i$ of an unknowm constant mean $m$ of
the field (see~Fig.~\ref{Fig2}).

\section{Summary}
\noindent
It was shown that the (estimation) statistics of
the field $V_j; j~\subset i=1,\ldots,n$ with
an unknown constant mean $m$ and variance $\sigma^2$
$$
\hat{V}_j=\sum_{i=1}^n \omega^i_j V_i = \omega^i_j V_i
$$
that assumes -- the unbiasedness constraint~(\ref{uc})
$$
E\{V_j\}-E\{\omega^i_j V_i\}
=
0
$$
that assumes -- the minimization constraint~(\ref{mc})
$$
\frac{\partial D^2\{V_j-\omega^i_j V_j\}}{\partial \omega^i_j}
=
0
$$
given by the kriging system of equations~(\ref{ke})
$$
\begin{array}{cccccl}
{\underbrace{
\left[
\begin{array}{cccc}
\rho_{\it 11} & \ldots & \rho_{\it 1n} & 1 \\
\vdots & \ddots & \vdots & \vdots \\
\rho_{\it n1} & \ldots & \rho_{\it nn} & 1 \\
1 & \ldots & 1 & 0 \\
\end{array}
\right]}_{(n+1) \times (n+1)}}
&
\cdot
&
\underbrace{
\left[
\begin{array}{c}
\omega_j^{\it 1} \\
\vdots \\
\omega_j^{\it n} \\
\mu^{\it 1}_j \\
\end{array}
\right]
}_{(n+1) \times 1}
&
=
&
\underbrace{
\left[
\begin{array}{c}
\rho_{\it 1j} \\
\vdots \\
\rho_{\it nj} \\
1 \\
\end{array}
\right]
}_{(n+1) \times 1} 
\end{array}
$$
is the best linear unbiased generalized (estimation) statistics of
random field $V_j$
with minimized variance of the statistics
\begin{equation}
E\{[\omega^i_j V_i-m]^2\}
=
\mp\sigma^2
\left( \omega^i_j \rho_{ij}- \mu^{\it 1}_j \right)
\label{o2}
\end{equation}
and (minimized) 
\begin{equation}
E\{[V_j-\omega^i_j V_i]^2\}
=
\sigma^2
\left( 1 \pm \left(\omega^i_j \rho_{ij} + \mu^{\it 1}_j \right)\right)
\label{oo2}
\end{equation}
with the asymptotic property~(Appendix~\ref{sec:B})
$$
\lim_{n \rightarrow \infty}
\lim_{j \rightarrow \infty}
E\{[\omega^i_j V_i - m]^2\}=0
$$
and
$$
\lim_{n \rightarrow \infty}
\lim_{j \rightarrow \infty}
E\{[V_j-\omega^i_j V_i]^2\}
=
\sigma^2
$$
constrained once again from~(\ref{oo2}) on computer -- 
the third constraint of spatial statistics
$$
E\{[V_j-\omega^i_j V_i]^2\}
=
\sigma^2
=
E\{[V_j-m]^2\}
$$
is the numerical best linear unbiased generalized statistics
for finite $n$ at finite $j$ of
an unknown constant mean $m=E\{V_j\}$ of the field $V_j$ with
the numerical generalized least-squares estimator $\omega^i_j v_i$ of
an unknown constant mean of the field and its asymptotic disjunction
for $j \rightarrow \infty$ the classic generalized least-squares
estimator $\lim_{j \rightarrow \infty} \omega^i_j v_i$ of 
an unknown constant mean of the field.

\appendix

\section{The sign of the terms}
\label{sec:A}
\noindent
If for correlation matrix $\rho_{il};~i,l=1,\ldots,n$
that consists of unit diagonal elements~(see~(\ref{cf}))
and non-positive off-diagonal elements
holds
$$
\omega^i_j
\rho_{il}
\omega^l_j < 0
$$
like at $j \ge n+1$ for vector
$\rho_{ij};~i=1,\ldots,n$
that consists of non-positive correlations
holds
$$
\omega^i_j
\rho_{ij}
< 0
$$
then~(\ref{D})
$$
D^2\{R_j\}
=\sigma^2
+2 \sigma^2
\omega^i_j
\rho_{ij}
-\sigma^2
\omega^i_j
\rho_{il}
\omega^l_j
$$
for non-negative correlation function
$$
D^2\{R_j\}=
\sigma^2
-2 \sigma^2 \omega^i_j \rho_{ij}
+ \sigma^2 \omega^i_j \rho_{il} \omega^l_j
$$
for white noise
$$
D^2\{R_j\}=
\sigma^2
+\sigma^2 \omega^i_j \rho_{il} \omega^l_j \ ,
$$
where $\rho_{il}$ is the identity matrix.

\section{The asymptotic property of the best
linear unbiased generalized statistics of random field}
\label{sec:B}
\noindent
From the minimization constraint~(\ref{mc})
$$
\rho_{is} \omega^s_j
+f_{\it i1} \mu^{\it 1}_j
=
\rho_{ij}
$$
we get
$$
\delta^l_s \omega^s_j
=
\omega^l_j
=
-\rho^{li} f_{\it i1} \mu^{\it 1}_j
+
\rho^{li} \rho_{ij} \ ,
$$
where
$$
\rho^{li}
\rho_{is}
=\delta^l_s \ ,
$$
substituted into~(\ref{uc})
$$
f_{\it 1l} \omega^l_j=1
$$
gives
$$
\mu^{\it 1}_j
=
(f_{\it 1l} \rho^{li} f_{\it i1})^{-1}
(f_{\it 1l} \rho^{li} \rho_{ij}-1)
$$
for $j \rightarrow \infty$~(\ref{cv})
$$
\rho_{ij}= \xi f_{\it i1}
$$
then the Lagrange parameter simplifies to
$$
\lim_{j \rightarrow \infty} \mu^{\it 1}_j
=\xi - (f_{\it 1l}\rho^{li} f_{\it i1})^{-1}
$$
from the unbiasedness condition~(\ref{uc})
it also holds
$$
\lim_{j \rightarrow \infty} \omega^i_j \rho_{ij}=
\xi \omega^i_j f_{\it i1} = \xi
$$
then the minimized variance of the estimation statistics~(\ref{o})
$$
E\{[\omega^i_j V_i-m]^2\}
=
\mp\sigma^2(\omega^i_j \rho_{ij}-\mu^{\it 1}_j)
$$
simplifies to
$$
\lim_{j \rightarrow \infty} E\{[\omega^i_j V_i-m]^2\}
=
\mp\sigma^2\left(\xi - \xi +
(f_{\it 1l} \rho^{li} f_{\it i1})^{-1} \right)
$$
and~(\ref{oo})
$$
E\{[V_j-\omega^i_j V_i]^2\}
=
\sigma^2
(1 \pm (\omega^i_j \rho_{ij}+\mu^{\it 1}_j))
$$
simplifies to
$$
\lim_{j \rightarrow \infty}E\{[V_j-\omega^i_j V_i]^2\}
=
\sigma^2
(1 \pm (\xi + \xi -
(f_{\it 1l} \rho^{li} f_{\it i1})^{-1}))
$$
since
$$
\lim_{n \rightarrow \infty}
\left(f_{\it 1l} \rho^{li} f_{\it i1} \right)^{-1}
=
0
$$
and
$$
\xi \rightarrow 0
$$
we get the asymptotic property of the best linear
unbiased generalized statistics of random field
$$
\lim_{n \rightarrow \infty}
\lim_{j \rightarrow \infty}
E\{[\omega^i_j V_i-m]^2\}
=
0
$$
and
$$
\lim_{n \rightarrow \infty}
\lim_{j \rightarrow \infty}
E\{[V_j-\omega^i_j V_i]^2\}
=
\sigma^2 \ .
$$

\end{document}